# Deconvolution density estimation with heteroscedastic errors using SIMEX

**Xiao-Feng Wang**[*][†]

*Department of Quantitative Health Sciences,
Cleveland Clinic, Cleveland, OH 44195, USA*
*e-mail:* `wangx6@ccf.org`

**Jiayang Sun**[*]

*Department of Statistics,
Case Western Reserve University, Cleveland, OH 44106, USA*
*e-mail:* `jiayang@sun.case.edu`

**Zhaozhi Fan**

*Department of Mathematics and Statistics,
Memorial University of Newfoundland, St. John's, NL A1C 5S7, Canada*
*e-mail:* `zhaozhi@math.mun.ca`

**Abstract:** In many real applications, the distribution of measurement error could vary with each subject or even with each observation so the errors are heteroscedastic. In this paper, we propose a fast algorithm using a simulation-extrapolation (SIMEX) method to recover the unknown density in the case of heteroscedastic contamination. We show the consistency of the estimator and obtain its asymptotic variance and then address the practical selection of the smoothing parameter. We demonstrate that, through a finite sample simulation study, the proposed method performs better than the Fourier-type deconvolution method in terms of integrated squared error criterion. Finally, a real data application is conducted to illustrate the use of the method.

## Contents

[*]The corresponding authors
[†]The research of Xiao-Feng Wang is supported in part by the NIH Awards: UL1 RR024989-01.

## 1. Introduction

A fundamental problem in measurement error models is to recover an unknown density of a variable when its observed values or data are contaminated with errors. The ever increasing interest in the problem comes from an increased number of medical and financial studies in which variables are observed with measurement errors or are only partial available. Formally, let $X$ be the variable of interest, which we cannot observe directly. Instead, based on an observed sample $Y_1, \cdots, Y_n$ drawn independently from

$$Y = X + U \tag{1}$$

where the measurement error $U$ is independent of $X$. One is interested in estimating the unknown density function of $X$. The distribution of $U$ is typically assumed known or can be estimated separately.

Denote the density functions of $X, U$ and $Y$ by $f_X, f_U$ and $f_Y$, and their characteristic functions by $\varphi_X, \varphi_U$ and $\varphi_Y$, respectively. If the measurement errors are ignored, then a *naive* estimator of $f_X(x)$ is the ordinary kernel estimate,

$$\hat{f}_{X,naive}(x) = \frac{1}{nh} \sum_{j=1}^{n} K\left(\frac{x-Y_j}{h}\right) = \frac{1}{n} \sum_{j=1}^{n} K_h(x - Y_j), \tag{2}$$

where $K_h(x) = K(x/h)/h$, $K(\cdot)$ is a symmetric probability kernel with a finite variance $\int x^2 K(x) < \infty$. It is clear that $\hat{f}_{X,naive}(x)$ is a biased estimator of $f_X(x)$ with

$$\mathrm{E}(\hat{f}_{X,naive}(x)) = f_X * f_U(x),$$

where "$*$" denotes convolution. Hence, finding a consistent estimator of $f_X(x)$ requires deconvoluting the density of measurement errors $f_U(x)$.

The usual *deconvolution* procedure is to apply a Fourier inversion on $\varphi_X(t)$:

$$f_X(x) = \frac{1}{2\pi} \int e^{-itx} \varphi_X(t) dt = \frac{1}{2\pi} \int e^{-itx} \frac{\varphi_Y(t)}{\varphi_U(t)} dt. \tag{3}$$

Of course, $\varphi_Y$ in (3), if unknown, needs to be estimated. Indeed, substituting $\varphi_Y(t)$ in (3) by its kernel estimate

$$\hat{\varphi}_Y(t) = \int e^{itx} \hat{f}_Y(x) dx = \frac{1}{nh} \sum_{j=1}^{n} \int e^{it(x-Y_i)} K\left(\frac{x-Y_j}{h}\right) e^{itY_j} dx = \frac{1}{n} \sum_{j=1}^{n} \varphi_K(th) e^{itY_j}$$

leads to the *Deconvoluting Kernel Estimator* (DKE) of $f_X$, first introduced by [21]:

$$\hat{f}_{X,DKE}(x) = \frac{1}{nh} \sum_{j=1}^{n} K^*\left(\frac{x-Y_j}{h}\right), \tag{4}$$

where

$$K^*(z) = \frac{1}{2\pi} \int e^{-itz} \frac{\varphi_K(t)}{\varphi_U(t/h)} dt \tag{5}$$

is called the *deconvoluting kernel*. Observe that the deconvoluting kernel estimate in (4) is just an ordinary kernel estimate but with a specific kernel function in (5). We call any estimator of $f_X(x)$ that requires a Fourier inversion or transformation a *Fourier-type estimator.*

There has been substantial literature on the Fourier-type deconvoluting kernel estimators. See, for instance, [21], [1], [15], [10], [9], [26], [5] and [25]. The difficulty with a deconvolution problem depends heavily on the smoothness of the error density $f_U$. Fan [10] characterized two types of error distributions: ordinary smooth and super-smooth distributions. Examples of the ordinary smooth distribution include gamma, symmetric gamma and Laplacian distributions; examples of the super-smooth distribution include normal, mixture normal and Cauchy distributions. The convergence rate of a DKE is very slow in the case of a super-smooth error distribution. For example, when $U$ is normally distributed, the convergence rate is only at $O\bigl((\log n)^{-1/2}\bigr)$ [27; 10].

Further, in many real applications, the distribution of measurement error could vary with each subject or even with each observation so the errors are heteroscedastic. Consideration of heteroscedastic errors can be traced back at least to Fuller [12] (chapter 3) in 1987, in a special case of linear regression, where predictor are measured with error. In a recent book, Carroll et al. [3] discussed systematically the state of art in measurement error models including the nonlinear regression with heteroscedastic error in variables. Delaigle and Meister[8] proposed an adjusted Fourier deconvolution estimator for the density estimation with heteroscedastic errors. They also applied the adjusted method to the nonparametric regression problem [7]. Staudenmayer et al. [20] considered a different type of model where the observed data are the sample means of replicates contaminated with heteroscedastic errors. They presented a spline-based density estimation method using a Monte Carlo Markov chain and a random-walk Metropolis-Hastings algorithm. Sun et al. [24] proposed new non-Fourier estimators for density function when errors are homoscedastic or heteroscedastic but uniformly distributed. The new estimators abandon the characteristic functions - there are no Fourier transformations needed in the calculation.

In this paper, we provide a new density estimation procedure for data contaminated with additive and heteroscedastic Gaussian measurement error (§2), without using a Fourier transform. Our resulting estimator is a "variable bandwidth type" kernel estimator, adopting the simulation-extrapolation (SIMEX) idea [22], though the simulation step in the original SIMEX algorithm is bypassed in our procedure. It is asymptotically unbiased and consistent (Proposition 2.1). The practical selection of the smoothing parameter is addressed in §3. Our estimator is computationally faster than the Fourier-type KDE proposed by [8]. It also has a competitive and often smaller integrated squared error than the DKE does (§4). An application to an astronomy data set using the both our method and the KDE method is given in §5. The article ends with concluding remarks in §6. The proof of Proposition 1 is in the Appendix.

## 2. Estimation procedure and asymptotics

SIMEX, first proposed by Cook and Stefanski [4], is a "jackknife"-type bias-adjusted method that has been widely applied in regression problems with measurement error. Cook and Stefanski [4; 23] applied the SIMEX algorithm to parametric regression problems. Carroll et al.[2]; Staudenmayer and Ruppert [19] discussed the nonparametric regression in the presence of measurement error using SIMEX method. Stefanski and Bay [22] focused on SIMEX estimation of a finite population cumulative distribution function (rather than a density) when sample units are measured with error. For a more complete discussion of the subject of SIMEX see the monograph by [3].

The key idea underlying SIMEX is the fact that the effect of measurement error on an estimator can be determined experimentally via simulation. SIMEX methods in both parametric and nonparametric regression consist of two steps: a simulation step and an extrapolation step. In the *simulation step*, additional independent measurement errors with known variance (typically denoted by $\lambda \sigma_U^2$, where $\lambda$ is a parameter to control the amount of added measurement errors) are generated and added to the original data, thereby creating "pseudo" data sets with successively larger measurement error variances. So, the total measurement error variance is then $(1 + \lambda_s)\sigma_U^2$ for the $s$th pseudo data set. The "pseudo estimators" are obtained from each of the generated "pseudo" data sets. The above simulation and estimation are repeated a large number of times, and the average value of the estimators for each level of contamination is calculated. Then, in the *extrapolation step*, a regression technique, such as, nonlinear least squares, is used to fit the trend between the pseudo estimators and the controlling parameter $\lambda$ of the added errors. At last, extrapolation to the ideal case of no measurement error (*i.e.* $\lambda = -1$) yields the SIMEX estimator.

To investigate a SIMEX algorithm on *density estimation* for data contaminated with errors, we consider a general *heteroscedastic measurement error model*. We generalize model (1) to

$$Y_j = X_j + U_j \qquad (6)$$

where $j = 1, \cdots, n$, $X_j \sim f_X$ and $U_j \sim f_{U_j}$. Since the normal distribution is frequently used in applications, we will focus on the typical super-smooth heteroscedastic error case: $U_j \sim N(0, \sigma_j^2)$, for $j = 1, \cdots, n$. Hence, $f_{U_j}$ is completely specified if $\sigma_j$ is. Clearly, if $\sigma_1 = \cdots = \sigma_n = \sigma$, errors are *homoscedastic*; otherwise, errors are *heteroscedastic*. Note that the homoscedastic model is just a special case of the heteroscedastic model.

Under the model setting (6), we now consider estimating the unknown density function $f_X(t)$ at a given $t$. By the general simulation extrapolation algorithm, estimators are re-computed on a large number $m$ of measurement error-inflated pseudo data sets, $\{\tilde{Y}_{jr}(\lambda)\}_{j=1}^n$, where

$$\tilde{Y}_{jr}(\lambda) = Y_j + \lambda^{1/2} U_{jr}, \qquad j = 1, \cdots, n, \quad r = 1, \cdots, m,$$

$U_{jr}$ are independent, pseudo-random variables with density $f_{U_j}$, and $\lambda \geqslant 0$ is a constant controlling the amount of added error.

The conventional kernel estimator of the density function from the $r$th variance-inflated data is

$$\hat{g}_r(t) = \frac{1}{nh} \sum_{j=1}^n K\left(\frac{t - \tilde{Y}_{jr}}{h}\right) = \frac{1}{n} \sum_{j=1}^n K_h(t - \tilde{Y}_{jr}).$$

By the SIMEX algorithm, the simulation and estimation steps are repeated a large number of times, and the average value of the estimators for each level of contamination is calculated.

$$\hat{g}(t) = \frac{1}{m} \sum_{r=1}^m \hat{g}_r(t) = \frac{1}{m} \sum_{r=1}^m \left(\frac{1}{n} \sum_{j=1}^n K_h(t - \tilde{Y}_{jr})\right) \qquad (7)$$

Note that by the law of large number,

$$\hat{g}(t) \xrightarrow{a.s.} \mathrm{E}(\hat{g}(t)) = \frac{1}{n} \sum_{j=1}^n \mathrm{E}\left(K_h(t - Y_j - \lambda^{1/2} U_{jr}) | Y_j\right)$$

as $m \to \infty$. Denote $\phi(\cdot)$ the density function of standard normal distribution and let $v = -(t - Y_j - \sigma_j \lambda^{1/2} w)/h$, we have as $h \to 0$,

$$\mathrm{E}\left(K_h(t - Y_j - \sigma_j \lambda^{1/2} W_{jr}) | Y_j\right)$$
$$= \int K_h(t - Y_j - \sigma_j \lambda^{1/2} w) \phi(w) dw$$
$$= \frac{1}{\sigma_j \lambda^{1/2}} \int K(v) \phi\left(\frac{t - Y_j + vh}{\sigma_j \lambda^{1/2}}\right) dv$$
$$= \frac{1}{\sigma_j \lambda^{1/2}} \int K(v) \left[\phi\left(\frac{t - Y_j}{\sigma_j \lambda^{1/2}}\right) + \frac{vh}{\sigma_j \lambda^{1/2}} \left(\frac{t - Y_j}{\sigma_j \lambda^{1/2}}\right) \phi\left(\frac{t - Y_j}{\sigma_j \lambda^{1/2}}\right) + o(h)\right] dv$$
$$= \frac{1}{\sigma_j \lambda^{1/2}} \phi\left(\frac{t - Y_j}{\sigma_j \lambda^{1/2}}\right) + o(h)$$

Therefore, the simulation step can be bypassed in our SIMEX algorithm for density estimation. The above derivation of (8) is similar to for the cumulative distribution function estimation by [22]. We use $g^*$ in (8) to replace $g$ in (7) in our estimation,

$$\hat{g}^*(t, \lambda) = \frac{1}{n} \sum_{j=1}^n \frac{1}{\sigma_j \lambda^{1/2}} \phi\left(\frac{t - Y_j}{\sigma_j \lambda^{1/2}}\right). \qquad (8)$$

We then calculate the quantity in (8) for a pre-determined sequence of $\lambda$, i.e. $0 \leq \lambda_1 < \lambda_2 < \cdots < \lambda_s$. The success of SIMEX technique depends on the assumption that the expectation of $\hat{g}^*$ is well-approximated by a nonlinear function of $\lambda$. Here we consider a standard quadratic function of $\lambda$,

$$E(\hat{g}^*(t,\lambda)) = \beta_0 + \beta_1 \lambda + \beta_2 \lambda^2.$$

Our SIMEX estimator for the unknown density $f$ without measurement error then can be obtained by the extrapolation step, *i.e.* letting $\lambda \to -1$

$$\hat{f}_{X,SIMEX}(t) = \lim_{\lambda \to -1} \widehat{E(\hat{g}^*(t,\lambda))} = \hat{\beta}_0 - \hat{\beta}_1 + \hat{\beta}_2. \tag{9}$$

Next, we can rewrite the SIMEX estimator in the form of

$$\hat{f}_{X,SIMEX}(t) = \frac{1}{n} \sum_{j=1}^{n} G_j(t|Y_j, \sigma_j, \Lambda), \tag{10}$$

where

$$G_j(t|Y_j, \sigma_j, \Lambda) = (1, -1, 1)(P^T P)^{-1} P^T \Phi_j(t|Y_j, \sigma_j, \Lambda),$$
$$P = (1, \Lambda, \Lambda^2), \ \Lambda = (\lambda_1, \cdots, \lambda_s)^T, \ \Lambda^2 = (\lambda_1^2, \cdots, \lambda_s^2)^T,$$

$$\Phi_j(t|Y_j, \sigma_j, \Lambda) = \left[ \frac{1}{\sigma_j \lambda_1^{1/2}} \phi\left(\frac{t-Y_j}{\sigma_j \lambda_1^{1/2}}\right), \cdots, \frac{1}{\sigma_j \lambda_s^{1/2}} \phi\left(\frac{t-Y_j}{\sigma_j \lambda_s^{1/2}}\right) \right]^T.$$

The following proposition shows that our SIMEX estimator (10) is an asymptotically unbiased estimator for the unknown density $f_X$.

**Proposition 2.1.** *Suppose that (a) the polynomial extrapolant is exact; (b) $f_X$ has a bounded and continuous derivative; and (c) $\sigma_j \leq \sigma_0$ for all $j$ and some $\sigma_0 < \infty$. Then the SIMEX estimator $\hat{f}_{X,SIMEX}(t)$ in (10) is a consistent and an asymptotically unbiased estimator of density function $f_X$ of unobserved sample $X$. The asymptotic variance of the estimator is*

$$Var(\hat{f}_{X,SIMEX}(t)) = \frac{f_X(t)}{n\sqrt{2\pi}\sigma_H} c(\Lambda) \Sigma_\Lambda c(\Lambda)^T,$$

*where $\sigma_H = \frac{n}{\sum_{j=1}^{n} \frac{1}{\sigma_j}}$, $c(\Lambda) = (1, -1, 1)(P^T P)^{-1} P^T$ and $\Sigma_\Lambda$ is a $s \times s$ matrix with the $lm^{th}$ element equals $\frac{1}{\sqrt{\lambda_l + \lambda_m}}$ ($l, m = 1, \cdots, s$).*

The proof of this proposition is given in the Appendix. □

*Remark* 2.1. It can be consistently estimated by replacing $f_X(t)$ with $\hat{f}_{X,SIMEX}(t)$ in the above formula. The asymptotic variance of $\hat{f}_{X,SIMEX}(t)$ equals the variance of a kernel estimator of $f(t)$ with adaptive bandwidth, multiplied by a constant which is not dependent on the random sample but the selected values of $\Lambda$. This fact enables us to choose an extrapolant that minimizes the asymptotic variance. The asymptotic normality of the SIMEX estimator is easily seen from the structure of the proposed SIMEX estimator. The confidence interval can be obtained correspondingly.

*Remark* 2.2. The asymptotics of the proposed estimator is based on the assumption that the polynomial extrapolant is exact. This is a typical assumption of SIMEX methods (see for example [2]), although it is difficult to verify it in real data analysis. We will show, in the simulation study, it is a reasonable even for complex models.

*Remark* 2.3. It is possible that the SIMEX density estimate takes negative values in small or sparse data regions, especially at tail regions. This behavior is common in many classes of kernel methods, such as wavelet density estimator, sinc kernel estimator, and spline estimator. This disadvantage will not affect the global performance of our estimate. A simple correction of our SIMEX estimator is

$$\hat{f}^*_{X,SIMEX}(t) = \max\{\hat{f}_{X,SIMEX}(t), 0\}.$$

## 3. Choice of the smoothing parameter

Our SIMEX estimator in (8) has the form of *variable-bandwidth kernel estimator* [14; 13], where $\sigma_j \lambda^{1/2}$ plays the role of the *variable-bandwidth* in the estimator. Thus, $\lambda$ can be considered as a smoothing parameter that determines the smoothness of the SIMEX estimating function. Our proposed SIMEX procedure requires specification of $0 \leq \lambda_1 < \lambda_2 < \cdots < \lambda_s$. The experience from our large simulation studies suggests that the choice of $s$ is not critical, neither is that of $\lambda_s$ if $\lambda_1$ is determined. A typical choice is taking $s = 50$, $\lambda_s = \lambda_1 + 3$ and $\lambda_2, \cdots, \lambda_{s-1}$ are equally spaced points in $(\lambda_1, \lambda_s)$. We propose two methods to select $\lambda_1$ here: one minimizes the mean integral squared error (MISE) and another is similar to the Silverman's rule-of-thumb method [18].

Based on Proposition 2.1, the MISE of the SIMEX density estimator, which depends on the parameter $\lambda$, is

$$\begin{aligned}
MISE(\hat{f}_{X,SIMEX}, \lambda) &= E(\int (\hat{f}_{X,SIMEX}(t) - f(t))^2 dt) \\
&= \int (Bias(\hat{f}_{X,SIMEX}(t)))^2 dt + \int Var(\hat{f}_{X,SIMEX}(t)) dt \\
&= \int Var(\hat{f}_{X,SIMEX}(t)) dt \approx \frac{c(\Lambda)\Sigma_\Lambda c(\Lambda)^T}{n\sqrt{2\pi}\sigma_H}.
\end{aligned}$$

Since only $\lambda_1$ is critical to the MISE of the SIMEX estimator, we propose to choose the parameter $\lambda_1$ and hence the bandwidth of the estimator through the minimization of $c(\Lambda)\Sigma_\Lambda c(\Lambda)^T$.

Another method to select the bandwidth is a type of rule-of-thumb method. For example, the bandwidth of a kernel density estimator is

$$\hat{h}_{Y,rot} = a_0 \min\left\{\hat{\sigma}_Y, \frac{R_Y}{1.34}\right\} n^{-1/5}$$

where $\hat{\sigma}_Y$ is the estimate of observed data $Y$, $R_Y = Y_{[0.75n]} - Y_{[0.25n]}$ is the inter-quartile range. Silverman's rule of thumb [18] uses factor $a_0 = 0.9$, while Scott [17] suggested using $a_0 = 1.06$. In the rest of the paper we use $a_0 = 1.06$.

From the equation (8), $\sigma_j \lambda_1^{1/2}$ plays the role of a *variable-bandwidth* in our estimator, based on the measurement error-inflated pseudo data sets $\tilde{Y}(\lambda_1)$. In order to determinate $\hat{\lambda}_{1,rot}$ we set:

$$\hat{h}_{\tilde{Y}(\lambda)} = \bar{\sigma}_U \hat{\lambda}_1^{1/2} = c_0 \hat{h}_Y$$

where we choose $c_0 = \sqrt{\sigma_Y^2 + \bar{\sigma}_U^2}/\sigma_Y$ and $\bar{\sigma}_U$ is the average of heteroscedastic standard deviation of measurement errors, *i.e.* $\bar{\sigma}_U = \sum_{j=1}^n \sigma_j/n$. So,

$$\hat{\lambda}_{1,rot} = \frac{\hat{\sigma}_Y^2 + \bar{\sigma}_U^2}{\hat{\sigma}_Y^2 \bar{\sigma}_U^2} \hat{h}_{Y,rot}. \tag{11}$$

## 4. Simulation study

In this section, the finite sample performance of the proposed SIMEX estimate is investigated via a simulation study. Our study involves the following three densities, representing typical features that can be encountered in practice:

(1) $X \sim N(0,1)$, a standard normal distribution.
(2) $X \sim Gamma(2,1)$, a Gamma distribution that is right skewed.
(3) $X \sim 0.5N(-2,1) + 0.5N(2,1)$, a normal mixture that is bimodal.

To assess the quality of a density estimator, we consider its *integrated squared error* (ISE) to the true density $f$:

$$ISE(\hat{f}) = \int \{\hat{f} - f(x)\}^2 dx.$$

The means and their standard errors of the ISEs of our SIMEX estimate, Fourier-type DKE, the naive estimate as well as the estimate obtained by using uncontaminated $X$ are compared below under the three typical distributions above for samples of sizes 50, 100, 250, and 1000, respectively. Figures of estimated curves below provide detailed pictures of the estimates in the entire domain of $X$. We have carried out more extensive simulation experiments by considering more complex target densities and other selections of the error variances than we can present here. Fortunately, all simulation experiments showed similar conclusions about the performance of the SIMEX estimator. All algorithms and simulations were implemented in R/Splus. Full results of the simulation study and R functions can be obtained from the authors upon request.

### 4.1. The case of homoscedastic errors

From each case of the densities, 1000 samples of size $n = 50, 100, 250,$ and 1000 were generated, each of which was then contaminated by a sample from a normal density of $N(0, \sigma_U^2)$. For each configuration, the parameter $\sigma_U$ was chosen equal to 0.2, 0.4, 0.6, 0.8, respectively.

*Choice of bandwidth.* From the studies of DKE, it is known that the kernel should be selected among densities whose characteristic function has a compact and symmetric support [11; 26; 5]. The use of such kernels guarantees the existence of the density estimator of (4). The most common example of such a kernel is given by the second-order kernel characteristic function

$$\phi_K(t) = (1-t^2)^3 I_{[-1,1]}(t),$$

where $I_{[-1,1]}(t)$ is the indicator function. The corresponding kernel function is

$$K(x) = \frac{48\cos x}{\pi x^4}\left(1 - \frac{15}{x^2}\right) - \frac{144\sin x}{\pi x^5}\left(2 - \frac{5}{x^2}\right).$$

This was the kernel we used in our simulation.

*Choice of bandwidth.* Since, Delaigle and Gijbels [6] showed that an appropriately chosen plug-in asymptotical bandwidth selector performed better than a cross-validated bandwidth selector and other bandwidth selectors, we therefore used the plug-in asymptotical bandwidth

$$\hat{h} = c_0 \sigma_W (\log n)^{-1/2},$$

following [11; 6], where $c_0 = 1.05$ in our simulation.

Results are shown in Table 1, in which entries without parentheses are the means and entries with parentheses are the standard errors of ISEs from our simulations of size 1000. The SIMEX estimate performs uniformly better than the DKE in terms of the ISE criterion. When the sample size increases, the means of ISEs of both the SIMEX estimate and the DKE get closer to the means of ISEs of kernel estimate of uncontaminated sample. The SIMEX estimate works beautifully for the cases of the moderate sample size and the large error variance. The standard errors of both the SIMEX estimate and the DKE are in the reasonable range.

[Table 1 about here.]

Figure 1 shows an example of deconvolution density estimation for the case of the homoscedastic measurement errors. $X$ is generated from $N(0,1)$ and two levels of errors $\sigma_U = 0.2$ and $\sigma_U = 0.8$ and four levels of sample size, $n = 50, 100, 250$, and $1000$ are considered in our study. In the figure, solid line denotes kernel estimate by uncontaminated sample X; dashed line denotes estimate by SIMEX method; dotted line denotes estimate by DKE method. Both the SIMEX and DKE methods recover the modes and capture the shape of true densities for the large sample sizes accurately. However, the DKE method is not stable when the sample size become smaller. With the small error variance (the sub-plot (a)), the DKE method shows wiggly curves for the small and modest sample sizes. This is due to the selection of the support kernel. A similar situation was also noted by [11]. A careful selection of the kernel function of DKE may improve the results. With the cases of the large error variance and/or small sample sizes, the DKE tends to underestimate the peaks while SIMEX works better.

[Figure 1 about here.]

### *4.2. The case of heteroscedastic errors*

In the case of heteroscedastic errors, measurement errors were generated from $N(0, \sigma_j^2)$, where $\sigma_j$ ($j = 1, \cdots, n$) were generated from a uniform distribution, $U(a, b)$. $(a, b)$ was chosen to be $(0.2, 0.4)$, $(0.4, 0.6)$, $(0.6, 0.8)$, $(0.8, 1)$, respectively. Due to the heavy computational burden of the Fourier-type method for the case of heteroscedastic errors, we considered a simulation study of size 500 rather than 1000 for each of the cases under sample sizes $n = 50, 100, 250$, and 1000, respectively. It took about several hours to finish a single case study using the Delaigle and Meister's Fourier-type method while our SIMEX procedure only took a few minutes to finish it. Of course, more efficient algorithm may be developed to speed up the computation of a Fourier-type estimate, but that was not our objective. Our simulation study of size 500 was already informative in comparing the performance of the SIMEX procedure and the Fourier-type procedure as shown in Table 2 and Figures 2 and 3.

The Fourier-type estimate we compared with was Delaigle and Meister's adjusted Fourier estimate for the density estimation with heteroscedastic error [8]. Their adjust estimator can be written as a form of a kernel density estimator,

$$\tilde{f}_{X,DKE}(x) = \frac{1}{nh} \sum_{j=1}^{n} \tilde{K}_j^* \left( \frac{x - Y_j}{h} \right),$$

where

$$\tilde{K}_j^*(z) = \frac{1}{2\pi} \int e^{-itz} \frac{\varphi_K(t)}{\psi_{U_j}(t/h)} dt, \quad \psi_{U_j}(t) = \frac{\frac{1}{n}\sum_{k=1}^{n} |\varphi_{U_k}(t)|^2}{\varphi_{U_j}(-t)}.$$

Table 2 shows the means and the standard errors of ISEs of the SIMEX estimate, Delaigle and Meister's Fourier-type estimate, the naive estimate and the standard kernel estimate based on uncontaminated data from 500 simulated samples for different cases. Similarly to the case of homoscedastic errors, the simulation results show that the SIMEX method performs uniformly better than the adjust DKE method in terms of the ISE. Comparing with Table 1, we find that the means and standard errors in the case of heteroscedastic errors are slightly larger than those in the case of homoscedastic errors. The results are not unexpected because heteroscedasticity of measurement errors brings more uncertainty and variation in the estimation than the homoscedastic errors.

[Table 2 about here.]

Figure 2 and Figure 3 display two examples of deconvolution density estimation in the case of the heteroscedastic errors. Two types of distributions are considered: a $Gamma(2, 1)$ distribution, and a normal mixture $0.5N(-2, 1) + 0.5N(2, 1)$. The errors are from $N(0, \sigma_j^2)$ where (a) $\sigma_j \sim U(0.2, 0.4)$, or (b) $\sigma_j \sim U(0.8, 1)$. The SIMEX estimates are closer than the true densities than the adjust DKE estimates do, especially at peaks and valleys. Despite the complex models, we see that the SIMEX estimate performs quite well in recovering the true densities and they are even competitive to the kernel estimate that are based on uncontaminated data.

[Figure 2 about here.]

[Figure 3 about here.]

## 5. A real data application

Our research was motivated with a real data analysis for astronomical data measured with errors. It is known that most astronomical data come with information subject to measurement errors. Sun et al. [24] gave an excellent introduction on the motivation of the measurement error problems in astronomy. Studying the distribution of one-dimensional velocities of stars originating in a given galaxy is of interest in astronomy. In this section, we illustrate our proposed method with an application to the astronomical position-velocity data in [16] from a sample of 26 low surfaces brightness (LSB) galaxies. The data contain 510 observed velocities of stars in km/s (relative to center, corrected for inclination) from 26 LSB galaxies. Each of observations includes its estimated standard deviation of measurement error. The sub-plot (a) of Figure 4 displays the histogram of measurement error standard deviations $\sigma_j$. The standard deviations vary from $0.1 \sim 46.8$ km/s and the mean is 6.34. The distribution is obviously skewed.

If we ignore the measurement errors, the velocity data look quite normal. The data range from $-289.00$ to $300.20$ and the mean is $-1.41$ and the median is $-1.00$. We applied the SIMEX method and the adjust DKE method to the data. The resulting estimated densities are shown in the sub-plot (b) of Figure 4. The two corrected estimates are consistent to each other, but not to the naive estimate. The probability around zero is higher and a small bump is detected on the left side of the curves (where the velocity approximately equal to $-250$ km/s) by both corrected estimates. This cannot be clearly seen from the naive estimate. Astronomers are mostly interested in this substructure of galaxies and can conduct further studies, which could lead to new discoveries.

[Figure 4 about here.]

## 6. Discussion

We presented a fast algorithm using SIMEX method to recover the unknown density when the data are contaminated with heteroscedastic errors and compared it with the Fourier-type method. The SIMEX estimate has advantages over the Fourier-type method in terms of ISE and computational efficiency and burden. We did not directly compare the SIMEX method with Bayesian method by Staudenmayer et al. [20], because they discussed the their method under a different model setting where the observed data are the sample means of replicates contaminated with heteroscedastic errors. We note that the Bayesian method has the similar simulation performance as the adjust DKE method in terms of ISE criterion under the model model setting of [20] (See table 1 of their paper), and it is also computational intensive.

Although the Fourier method has many nice theoretical properties, the SIMEX method is an excellent alternative in real data analysis. It is easy to implement and computationally more efficient. For example. The SIMEX method can also be used in classification analysis of microarray data with measurement error. In such a case, applying Fourier-type method is slow because we are facing deconvolution problems with thousands of genes.

## Appendix

*Proof of Proposition 2.1*:
  By (8),

$$
\begin{aligned}
E(\hat{g}^*(t,\lambda)) &= \frac{1}{n}\sum_{j=1}^{n}\frac{1}{\sigma_j\lambda^{1/2}}E\left[\phi(\frac{t-X_j-\sigma_j U_j}{\sigma_j\lambda^{1/2}})\right] \\
&= \frac{1}{n}\sum_{j=1}^{n}\frac{1}{\sigma_j\lambda^{1/2}}\iint_{R^2}\phi(\frac{t-x-\sigma_j u}{\sigma_j\lambda^{1/2}})f_X(x)\phi(u)dxdu \\
&= \frac{1}{n}\sum_{j=1}^{n}\iint_{R^2}f_X(t-\sigma_j u-\sigma_j\lambda^{1/2}v)\phi(v)\phi(u)dvdu \quad (12) \\
&= \frac{1}{n}\sum_{j=1}^{n}E\left[f_X(t-\sigma_j(1+\lambda)^{1/2}Z)\right],
\end{aligned}
$$

where $Z \sim N(0,1)$. Hence under conditions (b) and (c),

$$\lim_{\lambda\to-1}E(\hat{g}^*(t,\lambda)) = \lim_{\lambda\to-1}\frac{1}{n}\sum_{j=1}^{n}E\left[f_X(t-\sigma_j(1+\lambda)^{1/2}U)\right] = f_X(t).$$

So, by (9), $\hat{f}_{X,SIMEX}(t) = \hat{\beta}_0 - \hat{\beta}_1 + \hat{\beta}_2 \xrightarrow{pr} \beta_0 - \beta_1 + \beta_2 = f_X(t)$. Under the condition (a), $\hat{f}_{X,SIMEX}(t)$ is thus asymptotically unbiased and consistent.

The asymptotic variance can be calculated from the extrapolation step. Indeed, $\hat{g}^*(t,\lambda)$ can be treated as a kernel estimator of $g(t)$ with adaptive bandwidth $\sigma_j\lambda^{1/2}$ as shown in (8). Hence the asymptotic variance of this estimator is

$$
\begin{aligned}
Var(\hat{g}^*(t,\lambda)) &= \frac{1}{n^2}\sum_{j=1}^{n}\frac{1}{\sigma_j\lambda^{1/2}}E(\phi(V)f_X(t-\sigma_j(U+\lambda^{1/2}V))) \\
&\quad - \frac{1}{n^2}\sum_{j=1}^{n}\left[E(f_X(t-\sigma_j(1+\lambda)^{1/2}V))\right]^2,
\end{aligned}
$$

by tracing the same arguments as those in (12), where $U$ and $V$ are independent standard normal random variables. When the measurement error variances are small, the above variance can be approximated by

$$Var(\hat{g}^*(t)) = \frac{f_X(t)}{n\sigma_H \lambda^{1/2}} \frac{1}{\sqrt{2\pi}} + O(\frac{\bar{\sigma}_U}{n}),$$

where $\bar{\sigma}_U = \frac{1}{n}\sum_{j=1}^n \sigma_j$ is the average of the standard deviations, and $\sigma_H = n/\sum_{j=1}^n \frac{1}{\sigma_j}$ is the harmonic average of the standard deviations. The harmonic average, from its definition, characterizes the underlying data in a way similar to their minimum. The approximation error of the above calculation is $O(\frac{\bar{\sigma}_U}{n})$.

The covariance between $\hat{g}^*(t, \lambda)$, when $\lambda$ takes different values $\lambda_l$ and $\lambda_m$, is

$$cov(\hat{g}^*(t, \lambda_l), \hat{g}^*(t, \lambda_m))$$
$$= \frac{1}{n^2} \sum_{j=1}^n \frac{1}{\sqrt{2\pi}\sigma_j\sqrt{\lambda_l + \lambda_m}} E\left[f_X(t - \sigma_j(1 + \frac{\lambda_l \lambda_m}{\lambda_l + \lambda_m})^{1/2}W)\right]$$
$$- \frac{1}{n^2} E\left[f_X(t - \sigma_j(1 + \lambda_l)^{1/2}V)\right] E\left[f_X(t - \sigma_j(1 + \lambda_m)^{1/2}V)\right]$$
$$= \frac{f_X(t)}{n\sqrt{2\pi}\sigma_H} \frac{1}{\sqrt{\lambda_l + \lambda_m}} + O(\frac{\bar{\sigma}_U}{n}).$$

Let $c(\Lambda) = (1, -1, 1)(P^T P)^{-1} P^T$. Therefore, by (10), the asymptotic variance of our SIMEX density estimator is

$$Var(\hat{f}_{X,SIMEX}(t)) = \frac{1}{n^2} \sum_{j=1}^n c(\Lambda) cov(\Phi_j(t|Y_j, \sigma_j, \Lambda)) c(\Lambda)^T$$
$$= \frac{f_X(t)}{n\sqrt{2\pi}\sigma_H} c(\Lambda) \Sigma_\Lambda c(\Lambda)^T,$$

where $\Sigma_\Lambda$ be a $s \times s$ matrix with the $lm^{\text{th}}$ element equals $\frac{1}{\sqrt{\lambda_l + \lambda_m}}$. We complete the proof.

TABLE 1
*The means and standard errors of ISEs for simulation in the case of homoscedastic error. Simulation size is 1000. Entries without parentheses are the means and entries with parentheses are the standard errors of ISEs.*

| Density | $n$ | ISE | | | | | | | |
|---|---|---|---|---|---|---|---|---|---|
| | | SIMEX | fx | fy | DKE | SIMEX | fx | fy | DKE |
| Normal | | | $\sigma_U=0.2$ | | | | $\sigma_U=0.4$ | | |
| | 50 | 0.01040 | 0.01085 | 0.01147 | 0.02149 | 0.01670 | 0.01114 | 0.01333 | 0.01160 |
| | | (0.00030) | (0.00026) | (0.00027) | (0.00040) | (0.00044) | (0.00026) | (0.00027) | (0.00027) |
| | 100 | 0.00604 | 0.00647 | 0.00686 | 0.01331 | 0.00978 | 0.00646 | 0.00885 | 0.00784 |
| | | (0.00017) | (0.00015) | (0.00016) | (0.00024) | (0.00023) | (0.00013) | (0.00017) | (0.00018) |
| | 250 | 0.00300 | 0.00327 | 0.00361 | 0.00677 | 0.00483 | 0.00308 | 0.00520 | 0.00442 |
| | | (0.00007) | (0.00006) | (0.00007) | (0.00011) | (0.00011) | (0.00006) | (0.00009) | (0.00009) |
| | 1000 | 0.00103 | 0.00108 | 0.00135 | 0.00234 | 0.00172 | 0.00111 | 0.00289 | 0.00215 |
| | | (0.00002) | (0.00002) | (0.00002) | (0.00003) | (0.00004) | (0.00002) | (0.00004) | (0.00004) |
| | | | $\sigma_U=0.6$ | | | | $\sigma_U=0.8$ | | |
| | 50 | 0.02159 | 0.01148 | 0.01776 | 0.01638 | 0.02438 | 0.01151 | 0.02446 | 0.02902 |
| | | (0.00059) | (0.00028) | (0.00034) | (0.00030) | (0.00058) | (0.00026) | (0.00042) | (0.00033) |
| | 100 | 0.01230 | 0.00637 | 0.01312 | 0.01228 | 0.01422 | 0.00634 | 0.02033 | 0.02363 |
| | | (0.00030) | (0.00014) | (0.00024) | (0.00022) | (0.00031) | (0.00014) | (0.00029) | (0.00024) |
| | 250 | 0.00600 | 0.00314 | 0.00938 | 0.00837 | 0.00835 | 0.00329 | 0.01711 | 0.01817 |
| | | (0.00013) | (0.00006) | (0.00013) | (0.00012) | (0.00019) | (0.00006) | (0.00020) | (0.00017) |
| | 1000 | 0.00236 | 0.00112 | 0.00703 | 0.00547 | 0.00393 | 0.00113 | 0.01423 | 0.01270 |
| | | (0.00005) | (0.00002) | (0.00006) | (0.00006) | (0.00007) | (0.00002) | (0.00009) | (0.00008) |
| Gamma | | | $\sigma_U=0.2$ | | | | $\sigma_U=0.4$ | | |
| | 50 | 0.01222 | 0.01187 | 0.01250 | 0.02295 | 0.01526 | 0.01157 | 0.01469 | 0.01357 |
| | | (0.00030) | (0.00027) | (0.00027) | (0.00038) | (0.00037) | (0.00025) | (0.00028) | (0.00028) |
| | 100 | 0.00738 | 0.00716 | 0.00793 | 0.01388 | 0.00917 | 0.00698 | 0.01019 | 0.00906 |
| | | (0.00015) | (0.00014) | (0.00015) | (0.00021) | (0.00019) | (0.00013) | (0.00017) | (0.00016) |
| | 250 | 0.00420 | 0.00406 | 0.00477 | 0.00704 | 0.00551 | 0.00421 | 0.00737 | 0.00614 |
| | | (0.00007) | (0.00007) | (0.00007) | (0.00010) | (0.00011) | (0.00008) | (0.00010) | (0.00010) |
| | 1000 | 0.00181 | 0.00164 | 0.00228 | 0.00265 | 0.00255 | 0.00164 | 0.00451 | 0.00334 |
| | | (0.00002) | (0.00002) | (0.00003) | (0.00003) | (0.00004) | (0.00002) | (0.00004) | (0.00004) |
| | | | $\sigma_U=0.6$ | | | | $\sigma_U=0.8$ | | |
| | 50 | 0.01977 | 0.01143 | 0.01834 | 0.01676 | 0.02313 | 0.01195 | 0.02350 | 0.02512 |
| | | (0.00053) | (0.00025) | (0.00033) | (0.00028) | (0.00055) | (0.00027) | (0.00038) | (0.00031) |
| | 100 | 0.01234 | 0.00713 | 0.01426 | 0.01316 | 0.01589 | 0.00741 | 0.02006 | 0.02112 |
| | | (0.00028) | (0.00014) | (0.00022) | (0.00020) | (0.00033) | (0.00014) | (0.00026) | (0.00022) |
| | 250 | 0.00756 | 0.00406 | 0.01124 | 0.01008 | 0.01039 | 0.00400 | 0.01681 | 0.01689 |
| | | (0.00013) | (0.00007) | (0.00013) | (0.00012) | (0.00017) | (0.00007) | (0.00016) | (0.00014) |
| | 1000 | 0.00409 | 0.00164 | 0.00853 | 0.00705 | 0.00670 | 0.00166 | 0.01424 | 0.01298 |
| | | (0.00005) | (0.00002) | (0.00006) | (0.00005) | (0.00009) | (0.00003) | (0.00008) | (0.00007) |
| Mixture | | | $\sigma_U=0.2$ | | | | $\sigma_U=0.4$ | | |
| | 50 | 0.01448 | 0.01001 | 0.01058 | 0.02146 | 0.00864 | 0.01021 | 0.01224 | 0.01005 |
| | | (0.00011) | (0.00012) | (0.00012) | (0.00029) | (0.00015) | (0.00012) | (0.00013) | (0.00018) |
| | 100 | 0.00979 | 0.00715 | 0.00771 | 0.01299 | 0.00551 | 0.00727 | 0.00954 | 0.00675 |
| | | (0.00008) | (0.00008) | (0.00009) | (0.00016) | (0.00010) | (0.00008) | (0.00010) | (0.00012) |
| | 250 | 0.00472 | 0.00418 | 0.00473 | 0.00659 | 0.00264 | 0.00423 | 0.00643 | 0.00375 |
| | | (0.00005) | (0.00005) | (0.00005) | (0.00007) | (0.00005) | (0.00005) | (0.00006) | (0.00006) |
| | 1000 | 0.00115 | 0.00174 | 0.00216 | 0.00234 | 0.00088 | 0.00173 | 0.00357 | 0.00162 |
| | | (0.00002) | (0.00002) | (0.00002) | (0.00002) | (0.00002) | (0.00002) | (0.00003) | (0.00002) |
| | | | $\sigma_U=0.6$ | | | | $\sigma_U=0.8$ | | |
| | 50 | 0.01020 | 0.01033 | 0.01508 | 0.01103 | 0.01390 | 0.01033 | 0.01834 | 0.01664 |
| | | (0.00020) | (0.00012) | (0.00015) | (0.00019) | (0.00026) | (0.00012) | (0.00016) | (0.00017) |
| | 100 | 0.00651 | 0.00717 | 0.01218 | 0.00794 | 0.01005 | 0.00748 | 0.01609 | 0.01378 |
| | | (0.00013) | (0.00008) | (0.00011) | (0.00013) | (0.00019) | (0.00008) | (0.00012) | (0.00014) |
| | 250 | 0.00368 | 0.00428 | 0.00921 | 0.00541 | 0.00604 | 0.00431 | 0.01299 | 0.01029 |
| | | (0.00007) | (0.00005) | (0.00007) | (0.00008) | (0.00011) | (0.00005) | (0.00009) | (0.00010) |
| | 1000 | 0.00153 | 0.00172 | 0.00622 | 0.00324 | 0.00325 | 0.00175 | 0.01001 | 0.00713 |
| | | (0.00003) | (0.00002) | (0.00004) | (0.00003) | (0.00005) | (0.00002) | (0.00005) | (0.00005) |

Table 2
The means and standard errors of ISEs for simulation in the case of heteroscedastic error. Simulation size is 500. Entries without parentheses are the means and entries with parentheses are the standard errors of ISEs.

| Density | $n$ | ISE | | | | | | | |
|---|---|---|---|---|---|---|---|---|---|
| | | SIMEX | fx | fy | DKE | SIMEX | fx | fy | DKE |
| Normal | | $\sigma_U \sim U(0.2, 0.4)$ | | | | $\sigma_U \sim U(0.4, 0.6)$ | | | |
| | 50 | 0.01496 | 0.01110 | 0.01242 | 0.01422 | 0.01892 | 0.01093 | 0.01506 | 0.01275 |
| | | (0.00064) | (0.00038) | (0.00045) | (0.00047) | (0.00069) | (0.00037) | (0.00045) | (0.00042) |
| | 100 | 0.00846 | 0.00627 | 0.00736 | 0.00854 | 0.01173 | 0.00674 | 0.01076 | 0.00911 |
| | | (0.00027) | (0.00018) | (0.00022) | (0.00025) | (0.00039) | (0.00021) | (0.00029) | (0.00028) |
| | 250 | 0.00429 | 0.00323 | 0.00423 | 0.00452 | 0.00563 | 0.00320 | 0.00729 | 0.00597 |
| | | (0.00013) | (0.00009) | (0.00011) | (0.00012) | (0.00017) | (0.00009) | (0.00017) | (0.00016) |
| | 1000 | 0.00153 | 0.00111 | 0.00189 | 0.00170 | 0.00197 | 0.00113 | 0.00432 | 0.00306 |
| | | (0.00005) | (0.00003) | (0.00005) | (0.00005) | (0.00007) | (0.00004) | (0.00009) | (0.00008) |
| | | $\sigma_U \sim U(0.6, 0.8)$ | | | | $\sigma_U \sim U(0.8, 1)$ | | | |
| | 50 | 0.02318 | 0.01055 | 0.02039 | 0.02191 | 0.02687 | 0.01150 | 0.02995 | 0.03816 |
| | | (0.00101) | (0.00034) | (0.00058) | (0.00046) | (0.00090) | (0.00040) | (0.00066) | (0.00043) |
| | 100 | 0.01394 | 0.00666 | 0.01670 | 0.01735 | 0.01660 | 0.00629 | 0.02553 | 0.03167 |
| | | (0.00043) | (0.00021) | (0.00038) | (0.00034) | (0.00054) | (0.00018) | (0.00048) | (0.00036) |
| | 250 | 0.00685 | 0.00327 | 0.01269 | 0.01246 | 0.00977 | 0.00315 | 0.02153 | 0.02416 |
| | | (0.00020) | (0.00009) | (0.00023) | (0.00021) | (0.00028) | (0.00009) | (0.00015) | (0.00013) |
| | 1000 | 0.00285 | 0.00109 | 0.01008 | 0.00842 | 0.00552 | 0.00113 | 0.01885 | 0.01774 |
| | | (0.00009) | (0.00003) | (0.00013) | (0.00012) | (0.00016) | (0.00004) | (0.00019) | (0.00016) |
| Gamma | | $\sigma_U \sim U(0.2, 0.4)$ | | | | $\sigma_U \sim U(0.4, 0.6)$ | | | |
| | 50 | 0.01301 | 0.01128 | 0.01311 | 0.01571 | 0.01725 | 0.01113 | 0.01593 | 0.01397 |
| | | (0.00045) | (0.00034) | (0.00038) | (0.00042) | (0.00059) | (0.00031) | (0.00042) | (0.00039) |
| | 100 | 0.00861 | 0.00755 | 0.00942 | 0.00990 | 0.01072 | 0.00733 | 0.01247 | 0.01080 |
| | | (0.00025) | (0.00021) | (0.00024) | (0.00026) | (0.00031) | (0.00021) | (0.00027) | (0.00024) |
| | 250 | 0.00436 | 0.00391 | 0.00559 | 0.00521 | 0.00628 | 0.00386 | 0.00866 | 0.00723 |
| | | (0.00011) | (0.00010) | (0.00013) | (0.00012) | (0.00016) | (0.00010) | (0.00016) | (0.00015) |
| | 1000 | 0.00197 | 0.00162 | 0.00311 | 0.00230 | 0.00320 | 0.00165 | 0.00618 | 0.00474 |
| | | (0.00004) | (0.00003) | (0.00005) | (0.00004) | (0.00006) | (0.00003) | (0.00008) | (0.00007) |
| | | $\sigma_U \sim U(0.6, 0.8)$ | | | | $\sigma_U \sim U(0.8, 1)$ | | | |
| | 50 | 0.02210 | 0.01172 | 0.02163 | 0.02136 | 0.02471 | 0.01217 | 0.02686 | 0.03016 |
| | | (0.00077) | (0.00037) | (0.00051) | (0.00044) | (0.00080) | (0.00037) | (0.00055) | (0.00043) |
| | 100 | 0.01403 | 0.00735 | 0.01658 | 0.01621 | 0.01723 | 0.00739 | 0.02318 | 0.02578 |
| | | (0.00041) | (0.00021) | (0.00033) | (0.00029) | (0.00047) | (0.00020) | (0.00040) | (0.00032) |
| | 250 | 0.00896 | 0.00379 | 0.01391 | 0.01323 | 0.01196 | 0.00419 | 0.02018 | 0.02122 |
| | | (0.00023) | (0.00009) | (0.00021) | (0.00019) | (0.00030) | (0.00011) | (0.00026) | (0.00021) |
| | 1000 | 0.00518 | 0.00161 | 0.01098 | 0.00957 | 0.00796 | 0.00163 | 0.01735 | 0.01634 |
| | | (0.00010) | (0.00003) | (0.00011) | (0.00009) | (0.00014) | (0.00003) | (0.00013) | (0.00011) |
| Mixture | | $\sigma_U \sim U(0.2, 0.4)$ | | | | $\sigma_U \sim U(0.4, 0.6)$ | | | |
| | 50 | 0.00972 | 0.00994 | 0.01126 | 0.01259 | 0.00950 | 0.01042 | 0.01386 | 0.00999 |
| | | (0.00020) | (0.00016) | (0.00017) | (0.00029) | (0.00025) | (0.00017) | (0.00020) | (0.00027) |
| | 100 | 0.00644 | 0.00731 | 0.00870 | 0.00816 | 0.00592 | 0.00729 | 0.01086 | 0.00667 |
| | | (0.00013) | (0.00012) | (0.00013) | (0.00017) | (0.00015) | (0.00012) | (0.00015) | (0.00017) |
| | 250 | 0.00308 | 0.00416 | 0.00532 | 0.00430 | 0.00320 | 0.00424 | 0.00778 | 0.00415 |
| | | (0.00007) | (0.00006) | (0.00008) | (0.00008) | (0.00008) | (0.00007) | (0.00009) | (0.00009) |
| | 1000 | 0.00094 | 0.00175 | 0.00277 | 0.00153 | 0.00118 | 0.00170 | 0.00473 | 0.00210 |
| | | (0.00002) | (0.00003) | (0.00003) | (0.00003) | (0.000003) | (0.00003) | (0.00005) | (0.00004) |
| | | $\sigma_U \sim U(0.6, 0.8)$ | | | | $\sigma_U \sim U(0.8, 1)$ | | | |
| | 50 | 0.01140 | 0.01022 | 0.01641 | 0.01315 | 0.01581 | 0.01053 | 0.01997 | 0.01961 |
| | | (0.00031) | (0.00017) | (0.00022) | (0.00026) | (0.00043) | (0.00017) | (0.00024) | (0.00022) |
| | 100 | 0.00767 | 0.00716 | 0.01370 | 0.01000 | 0.01133 | 0.00719 | 0.01755 | 0.01685 |
| | | (0.00021) | (0.00011) | (0.00016) | (0.00019) | (0.00030) | (0.00012) | (0.00018) | (0.00017) |
| | 250 | 0.00498 | 0.00434 | 0.01133 | 0.00783 | 0.00776 | 0.00431 | 0.01508 | 0.01350 |
| | | (0.00012) | (0.00007) | (0.00012) | (0.00013) | (0.00017) | (0.00007) | (0.00013) | (0.00013) |
| | 1000 | 0.00117 | 0.00091 | 0.00412 | 0.00252 | 0.00454 | 0.00175 | 0.01202 | 0.00961 |
| | | (0.00006) | (0.00004) | (0.00018) | (0.00012) | (0.00008) | (0.00003) | (0.00007) | (0.00008) |

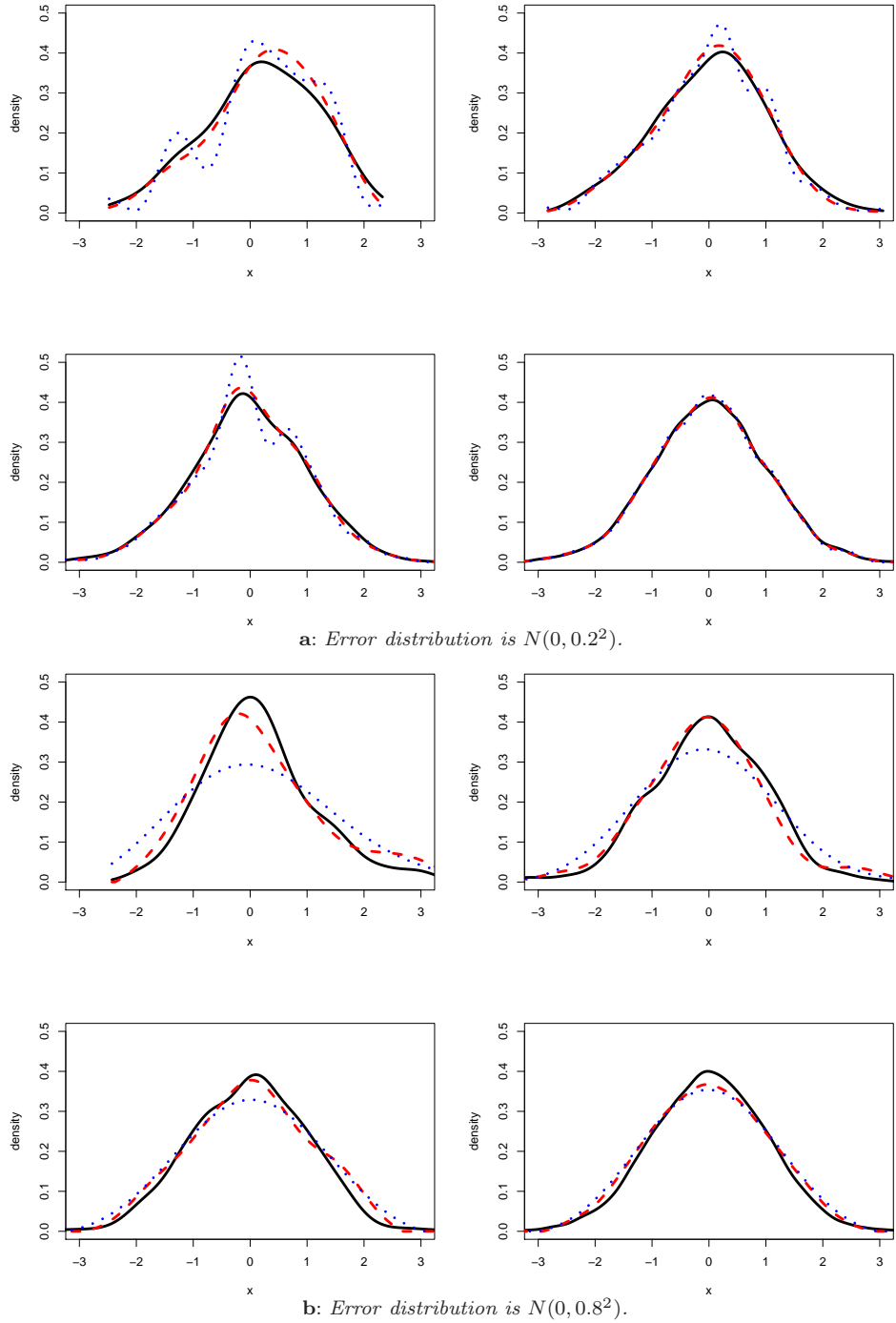

Fig 1. *Deconvolution density estimation in the case of heteroscedastic error: the true density is $N(0,1)$ and measurement errors are from (a) $N(0, 0.2^2)$ and (b) $N(0, 0.8^2)$ with different sample sizes. For both sub-plots (a) and (b), $n = 50$ (top left panel), $n = 100$ (top right panel), $n = 250$ (bottom left panel) and $n = 1000$ (bottom right panel). Solid line – kernel estimate by uncontaminated sample X; dashed line – estimate by SIMEX method; dotted line – estimate by DKE method.*

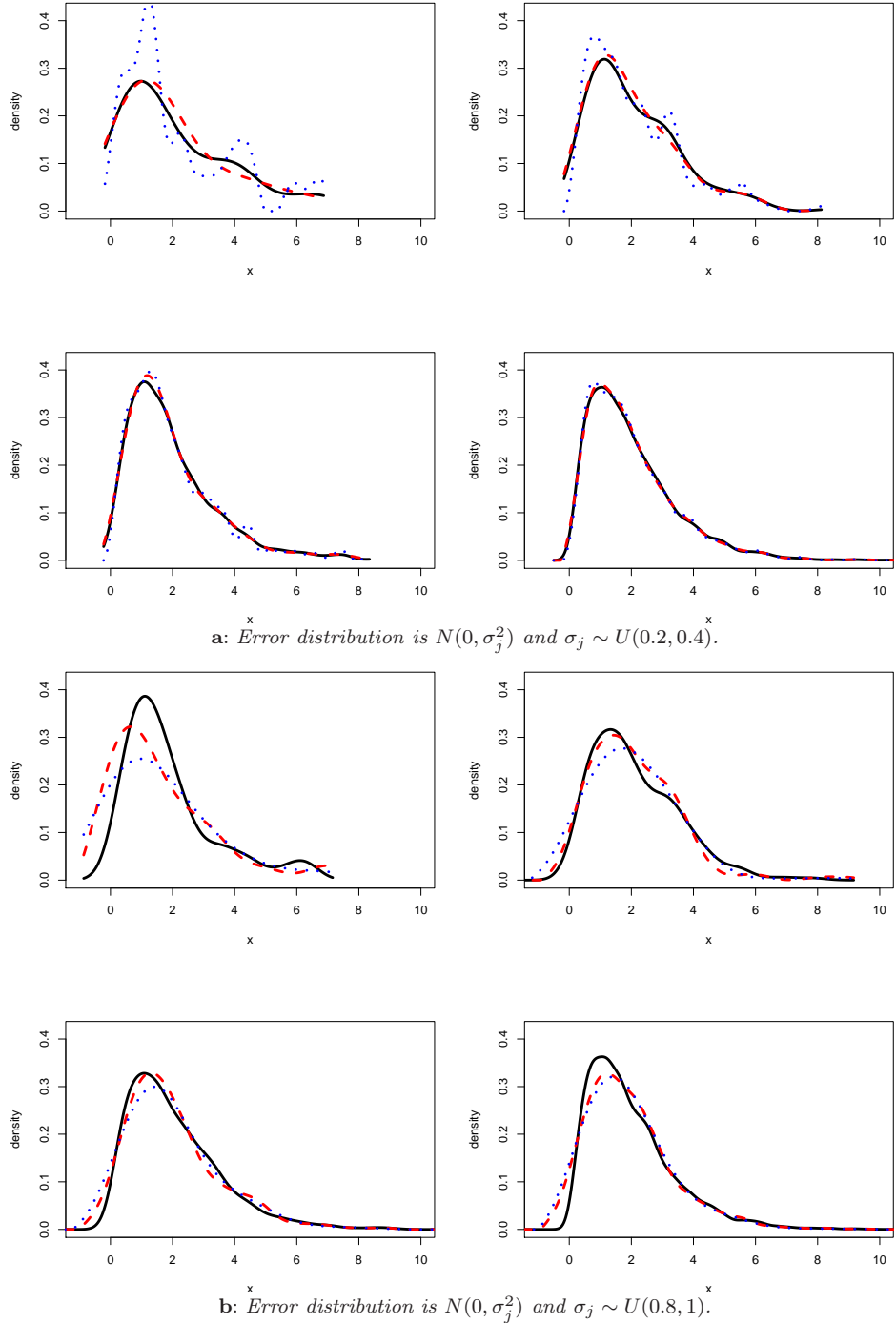

**a**: *Error distribution is $N(0, \sigma_j^2)$ and $\sigma_j \sim U(0.2, 0.4)$.*

**b**: *Error distribution is $N(0, \sigma_j^2)$ and $\sigma_j \sim U(0.8, 1)$.*

FIG 2. *Deconvolution density estimation in the case of heteroscedastic error: the true density is Gamma(2, 1) and measurement errors are from (a) $N(0, \sigma_j^2)$, $\sigma_j \sim U(0.2, 0.4)$ and (b) $N(0, \sigma_j^2)$, $\sigma_j \sim U(0.8, 1)$ with different sample sizes. For both subplots (a) and (b), $n = 50$ (top left panel), $n = 100$ (top right panel), $n = 250$ (bottom left panel) and $n = 1000$ (bottom right panel). Solid line – kernel estimate by uncontaminated sample X; dashed line – estimate by SIMEX method; dotted line – estimate by adjusted DKE method.*

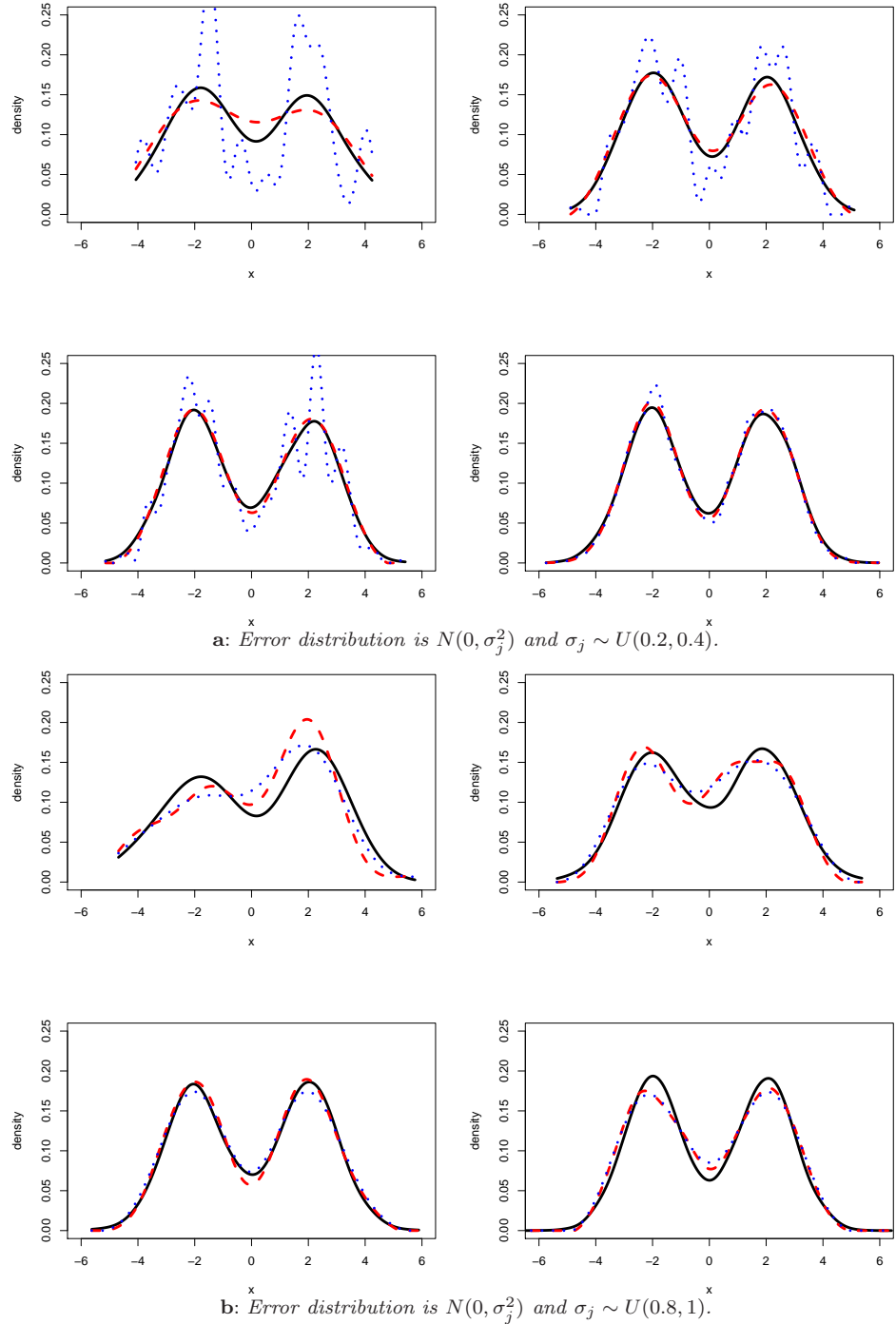

FIG 3. *Deconvolution density estimation in the case of heteroscedastic error: the true density is $0.5N(-2,1) + 0.5N(2,1)$ and measurement errors are from (a) $N(0, \sigma_j^2)$, $\sigma_j \sim U(0.2, 0.4)$ and (b) $N(0, \sigma_j^2)$, $\sigma_j \sim U(0.8, 1)$ with different sample sizes. For both sub-plots (a) and (b), $n = 50$ (top left panel), $n = 100$ (top right panel), $n = 250$ (bottom left panel) and $n = 1000$ (bottom right panel). Solid line – kernel estimate by uncontaminated sample X; dashed line – estimate by SIMEX method; dotted line – estimate by adjusted DKE method.*

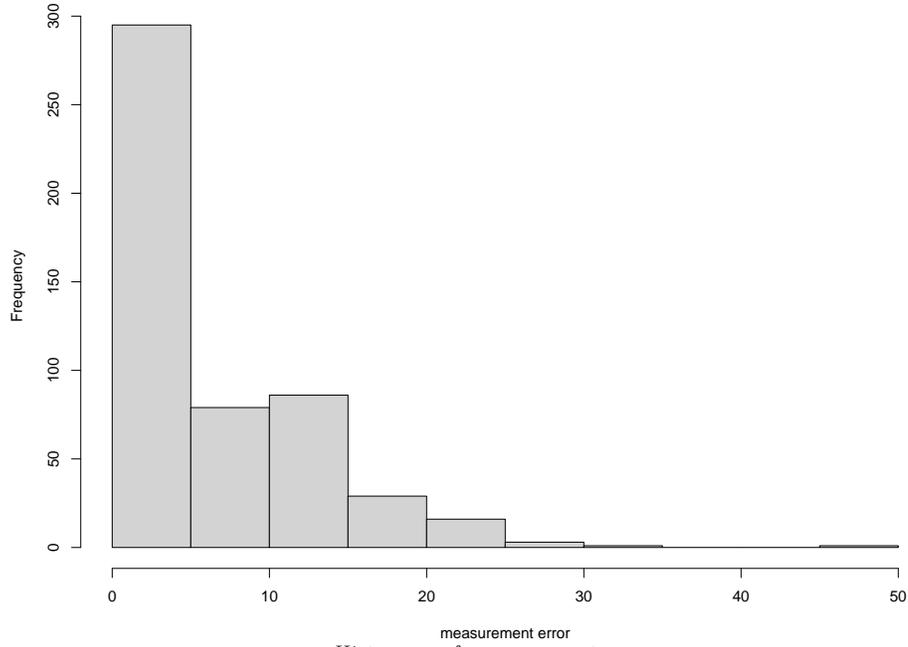

**a**: *Histogram of measurement errors.*

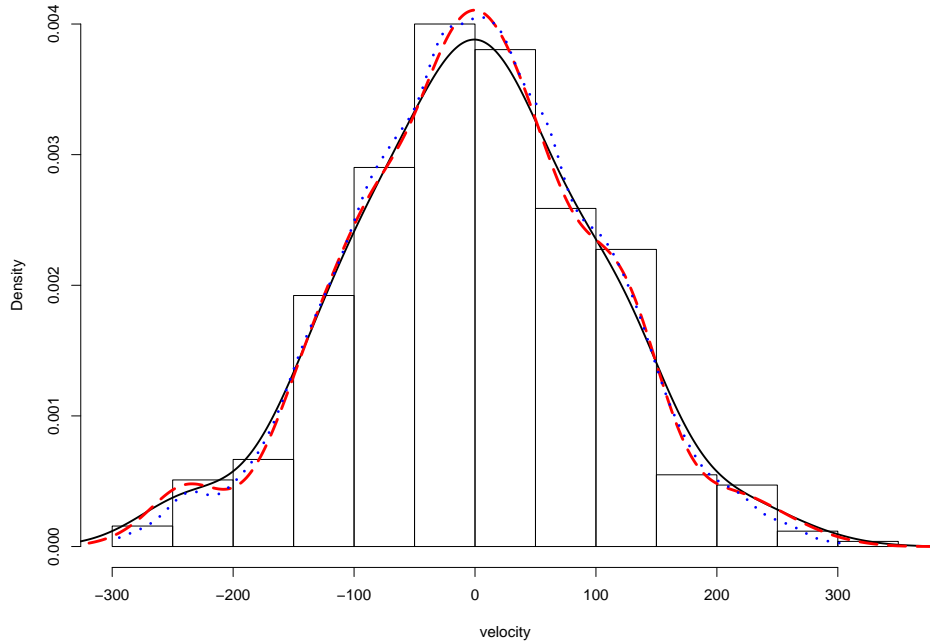

**b**: *Density estimation.*

Fig 4. *Density estimation of velocities in astronomical position-velocity data. The solid line is the naive estimate ignoring the heteroscedastic measurement errors. Two corrected estimates are considered here: SIMEX method (dashed line), DKE method (dotted line).*